\documentclass[11pt]{article}

\usepackage{amsmath,amssymb,amsthm,graphicx,enumerate}
\usepackage[letterpaper,margin=25mm]{geometry}
\usepackage[colorlinks,linkcolor=blue,urlcolor=blue,citecolor=blue]{hyperref}

\def\baselinestretch{1.25}

\newtheorem{theorem}{Theorem}
\newtheorem{lemma}{Lemma}

\newtheorem{problem}{Problem}
\newtheorem{conjecture}{Conjecture}

\def\etal{\emph{et al.}}

\def\le{\leqslant}
\def\ge{\geqslant}

\begin{document}
\title{On the Second-Order Wiener Ratios of Iterated Line Graphs}
\author{\textbf{Mohammad Ghebleh} and \textbf{Ali Kanso}\\[1ex]
Department of Mathematics\\ Kuwait University, P.O. Box 5969, Safat 13060, Kuwait\\[1ex]
{\sf\{mohammad.ghebleh, ali.kanso\}@ku.edu.kw}
}

\maketitle

\begin{abstract}
The Wiener index $W(G)$ of a graph $G$ is the sum of distances between all unordered pairs of its vertices.
Dobrynin and Mel'nikov~[in: Distance in Molecular Graphs -- Theory, 2012, p.\,\,85–121] propose the study of estimates for extremal values of the ratio $R_k(G)=W(L^k(G))/W(G)$ where $L^k(G)$ denotes the $k$th iterated line graph of~$G$.
Hri\v{n}\'akov\'a, Knor and \v{S}krekovski~[Art Discrete Appl.\,\,Math.\,\,1~(2018) \#P1.09] prove that for each $k\ge3$, the path $P_n$ has the smallest value of the ratio $R_k$ among all trees of large order~$n$, and they conjecture that the same holds for the case $k=2$. We give a counterexample of every order $n\ge22$ to this conjecture. 
\end{abstract}

\section{Introduction}

The Wiener index $W(G)$ of a graph $G$ is the sum of distances between all unordered pairs of its vertices. 
It was first introduced by Harry Wiener in 1947~\cite{wiener1947} as a structural descriptor of acyclic organic molecules.
Since the late 1970s, Wiener index has attracted the attention of graph theorists as a measure of the compactness of a graph. It is also referred to as \emph{distance of a graph}, \emph{average (mean) distance of a graph}, and \emph{transmission of a graph}~\cite{entringer1976,plesnik1984,soltes1991,chung1988,chung2002}.

The line graph $L(G)$ of a graph $G=(V,E)$ has $E$ as its vertex set, 
where $e,e'\in E$ are adjacent in $L(G)$ if and only if they have a common endvertex as edges of~$G$.
The iterated line graph $L^k(G)$ for a positive integer $k$ is defined via compositions of $L$ as a graph operator.
More specifically, $L^0(G)=G$, and $L^k(G)=L(L^{k-1}(G))$ for all $k\ge1$. 

The study of Wiener indices of iterated line graphs of a graph has raised considerable interest among graph theorists~\cite{dobr1,dobr2,dobr3,dou2010,edge-wiener,IterLine1,IterLine2,IterLine3,IterLine4,IterLine5,IterLine6,quipu1,quipu2}.
One of the earliest results of this type is due to Buckley~\cite{buckley1981} and states that for any tree $T$ of order $n$,
\[W(L(T))=W(T)-{n\choose2}.\]
Consequently, $W(L(T))<W(T)$ for all trees $T$ of order at least~$2$. On the other hand, for $k\ge2$, $W(L^k(T))$ may be smaller than, equal to, or larger than $W(T)$~\cite{IterLine1,IterLine2,IterLine3,IterLine4,IterLine5,IterLine6}.
Dobrynin and Mel'nikov~\cite{dobrynin2012} propose the study of estimates for extremal values of
the ratio \[R_k(G)=\frac{W(L^k(G))}{W(G)}.\]

The minimum value of $R_1$ is settled in Knor, \etal~\cite{tepeh2015}.
\begin{theorem}{\rm\cite{tepeh2015}}
Among all connected graphs $G$ of order~$n$, the ratio $R_1(G)$ is minimum for the star $S_n=K_{1,n-1}$.
\end{theorem}

Furthermore, for $k\ge3$, the minimum value of $R_k$ over trees of order $n$ is settled by Hri\v{n}\'akov\'a, \etal~\cite{hrinakova2018}.
\begin{theorem}{\rm\cite{hrinakova2018}}
Let $k\ge3$. Then the path $P_n$ attains the minimum value of $R_k$ in the class of trees on $n$ vertices.
\end{theorem}
Motivated by these results, the authors of~\cite{hrinakova2018} propose the following conjecture for the remaining cases.
\begin{conjecture}{\rm\cite{hrinakova2018}}
Let $n$ be a large number and $k\ge2$. Then among all graphs $G$ on $n$ vertices,
$R_k(G)$ attains its maximum at $G=K_n$, and its minimum at $G=P_n$.
\label{conj:main}
\end{conjecture}

We provide counterexamples of various homeomorphic classes for the minimum case of this conjecture when $k=2$.
That is, for every large enough integer $n$, we give a tree $T$ of order $n$ that satisfies
\begin{equation}
R_2(T)<R_2(P_n)
\label{eq:main}
\end{equation}
Our examples showcase the abundance of counterexamples to Conjecture~\ref{conj:main}, and
suggest that a complete characterization of trees (graphs) of a fixed order $n$ that minimize $R_2$ is likely to be much more complex than the cases~$k\neq2$.

\section{Background and Notation}

It is well-known that the Wiener index of the path of order $n$ equals 
\begin{equation}
W(P_n)=\tfrac16(n-1)n(n+1).
\label{eq:WPath}
\end{equation}
On the other hand, if $n\ge2$, then $L^2(P_n)=P_{n-2}$. Hence
\begin{equation}
R_2(P_n)=\frac{(n-2)(n-3)}{n(n+1)}=1-\frac{6(n-1)}{n(n+1)}.
\label{eq:Pn}
\end{equation}

We use the short-hand notation $W_k(G)=W(L^k(G))$ for the Wiener index of iterated line graphs of a graph~$G$.
In our study of $W_2(T)$ where $T$ is a tree, it is often convenient to compute the difference $D_2(T)=W(T)-W_2(T)$.
With this convention,
\[R_2(T)=\frac{W_2(T)}{W(T)}=\frac{W(T)-D_2(T)}{W(T)}=1-\frac{D_2(T)}{W(T)}.\]
In particular, to prove that a tree $T$ of order $n$ satisfies the inequality~(\ref{eq:main}), we may equivalently show
\begin{equation}
\frac{D_2(T)}{W(T)}>\frac{D_2(P_n)}{W(P_n)}=\frac{6(n-1)}{n(n+1)}.
\label{eq:mainD}
\end{equation}

\section{Counterexamples homeomorphic to the star $K_{1,3}$}

\begin{figure}
\begin{center}
\raisebox{-.5\height}{\includegraphics[scale=1.25,page=1]{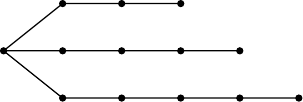}}
\hspace{3em}
\raisebox{-.5\height}{\includegraphics[scale=1.25,page=2]{figs.pdf}}
\end{center}
\caption{The tree $T_{3,4,5}$ and its second-order iterated line graph $L^2(T_{3,4,5})$.}
\label{fig:fan}
\end{figure}

For integers $a,b,c\ge1$, let $T_{a,b,c}$ denote the tree obtained from $K_{1,3}$ by subdividing (if necessary) its three edges to paths of length $a$, $b$ and $c$ respectively. Figure~\ref{fig:fan} shows $T_{3,4,5}$ and its second-order iterated line graph.

The Wiener index of these trees is obtained in~\cite{hrinakova2018}:
\begin{lemma}{\rm\cite{hrinakova2018}}
For positive integers $a,b,c$ we have \[W(T_{a,b,c})=\frac16(a+b+c)(a+b+c+1)(a+b+c+2) - abc.\]
\label{lemWTabc}
\end{lemma}

We will also need the Wiener index of the second-order iterated line graph of these trees.
\begin{lemma}
Let $a,b,c\ge2$ be integers. Then \[D_2(T_{a,b,c})=W(T_{a,b,c})-W_2(T_{a,b,c})=\frac12(a^2+b^2+c^2)+2(ab+ac+bc)-\frac12(a+b+c).\]
\label{lemW2Tabc}
\end{lemma}

\begin{proof}
Let $a,b,c\ge2$ be integers, and $T=T_{a,b,c}$.
We expand the Wiener index of $T$ and $L^2(T)$ as the sum of distances between pairs of their vertices to obtain
\begin{equation}
W(T)=W(P_{a+1})+W(P_{b+1})+W(P_{c+1})+\sum\limits_{i=1}^a\sum\limits_{j=1}^b(i+j)+\sum\limits_{i=1}^a\sum\limits_{j=1}^c(i+j)+\sum\limits_{i=1}^b\sum\limits_{j=1}^c(i+j),
\label{eq:WTabc}
\end{equation}
and
\begin{equation}
\begin{aligned}
W_2(T)=W(P_{a-1})+W(P_{b-1})+W(P_{c-1})&+\sum\limits_{i=1}^{a-1}(3i+1)+\sum\limits_{i=1}^{b-1}(3i+1)+\sum\limits_{i=1}^{c-1}(3i+1)+3\\
&+\sum\limits_{i=1}^{a-1}\sum\limits_{j=1}^{b-1}(i+j)+\sum\limits_{i=1}^{a-1}\sum\limits_{j=1}^{c-1}(i+j)+\sum\limits_{i=1}^{b-1}\sum\limits_{j=1}^{c-1}(i+j).
\end{aligned}
\label{eq:W2Tabc}
\end{equation}
Here for $W_2(T)$, we partition the vertices of $G=L^2(T)$ into the triangle on $\{x,y,z\}$ (as shown in Figure~\ref{fig:fan}), and the three connected components of $G-\{x,y,z\}$ (these are paths of lengths $a-2$, $b-2$, and $c-2$). It is immediate from equation~(\ref{eq:WPath}) that $W(P_{k+1})-W(P_{k-1})=k^2$ for all $k\ge2$.
On the other hand,
\[\sum\limits_{i=1}^{a-1}(3i+1)=a-1+\frac32a(a-1)=\frac32a^2-\frac12a-1,\]
and
\[\sum\limits_{i=1}^a\sum\limits_{j=1}^b(i+j)-\sum\limits_{i=1}^{a-1}\sum\limits_{j=1}^{b-1}(i+j)
=\sum\limits_{i=1}^{a}(i+b)+\sum\limits_{j=1}^{b}(a+j)-(a+b)=\frac12(a^2+b^2)+2ab-\frac12(a+b).\]
Similar calculations hold for the other sums involved in the equations~(\ref{eq:WTabc}) and~(\ref{eq:W2Tabc}).
Altogether these give
\begin{equation*}
D_2(T)=\tfrac12(a^2+b^2+c^2)+2(ab+ac+bc)-\tfrac12(a+b+c).
\qedhere
\end{equation*}
\end{proof}

As an example, for the tree $T=T_{7,7,7}$ of order $22$, Lemmas~\ref{lemWTabc} and~\ref{lemW2Tabc} give 
\[W(T)=1428 \text{ \ and \ } D_2(T)=357 ~\implies~ 1-R_2(T)=\frac14.\]
On the other hand, by equation~(\ref{eq:Pn}) we have
\[1-R_2(P_{22})=\frac{126}{506}<\frac14 \implies R_2(T)<R_2(P_n).\]
We conclude that $T=T_{7,7,7}$ is a counterexample (of order $22$) to Conjecture~\ref{conj:main}.
Indeed any \emph{balanced} tree $T_{a,b,c}$ where each two of the integers $a,b,c$ differ by at most~$1$ satisfies inequality~(\ref{eq:main}), provided that it has order at least~$22$. In the following theorem we prove this claim for large~$n$.

\begin{theorem}
There exists a number $n_0$ such that for all $n\ge n_0$, there exists a tree $T$ of order~$n$ satisfying $R_2(T)<R_2(P_n)$.
\label{thm:star}
\end{theorem}

\begin{proof}
Let $n\ge7$ be an integer and let $a=\lfloor(n-1)/3\rfloor$. Depending on the congruence class of $n$ modulo~$3$, we take $T$ to be one of the trees $T_{a,a,a}$, $T_{a,a,a+1}$, and $T_{a,a+1,a+1}$ that has order~$n$. 
\begin{enumerate}[(i)]
\item Suppose that $T=T_{a,a,a}$. Using Lemmas~\ref{lemWTabc} and~\ref{lemW2Tabc} we obtain
\[W(T)=\frac16(3a)(3a+1)(3a+2)-a^3=\frac12a(a+1)(7a+2),\]
and
\[D_2(T)=\frac12(3a^2)+2(3a^2)-\frac12(3a)=\frac32a(5a-1).\]
Therefore,
\[1=R_2(T)=\frac{3(5a-1)}{(a+1)(7a+2)}.\]
On the other hand, since in this case $n=3a+1$, equation~(\ref{eq:Pn}) gives
\[1-R_2(P_n)=\frac{18a}{(3a+1)(3a+2)}.\]
Now since
\[\lim_{a\to\infty}\frac{1-R_2(T)}{1-R_2(P_n)}=\lim_{a\to\infty}\frac{(3a+1)(3a+2)(5a-1)}{6a(a+1)(7a+2)}=\frac{15}{14},\]
there exists $n_1\in\mathbb{R}$ such that for all $n\ge n_1$ we have
\[\frac{1-R_2(T)}{1-R_2(P_n)}>1 ~\implies~ 1-R_2(T)>1-R_2(P_n) ~\implies~ R_2(T)<R_2(P_n).\]
\item If $T=T_{a,a,a+1}$, then $n=3a+2$ and similar calculations give
\[R_2(T)=1-\frac{a(15a+7)}{(a+1)^2(7a+2)} \text{ \ \ and \ \ } R_2(P_n)=1-\frac{2(3a+1)}{(a+1)(3a+2)},\]
and in turn,
\[\lim_{a\to\infty}\frac{1-R_2(T)}{1-R_2(P_n)}=\frac{15}{14},\]
which yields $R_2(T)<R_2(P_n)$ for all $n\ge n_2$ where $n_2\in\mathbb{R}$ is constant.
\item If $T=T_{a,a+1,a+1}$, then $n=3a+3$ and similar calculations give
\[R_2(T)=1-\frac{(3a+1)(5a+4)}{(a+1)(7a^2+16a+8)} \text{ \ \ and \ \ } R_2(P_n)=1-\frac{2(3a+2)}{(a+1)(3a+4)}.\]
and in turn,
\[\lim_{a\to\infty}\frac{1-R_2(T)}{1-R_2(P_n)}=\frac{15}{14},\]
which yields $R_2(T)<R_2(P_n)$ for all $n\ge n_3$ where $n_3\in\mathbb{R}$ is constant.
\end{enumerate}
We may take the largest of the constants $n_1,n_2,n_3$ to serve as the constant $n_0$.
\end{proof}

It is worth noting that as we have explicit formulae for $R_2(T)$ and $R_2(P_n)$ in all cases of the above proof,
a computer algebra system can be used to verify that in each case, the largest real root of the rational equation $R_2(T)=R_2(P_n)$ lies in the interval $(6,7)$. Thus we have $R_2(T)<R_2(P_n)$ in all cases provided that $a\ge7$.
This implies that in Theorem~\ref{thm:star}, the constant $n_0$ can be taken to be~$22$.

\section{Further homeomorphic classes}

In the preceding section we found trees of every large enough order $n$ that satisfy inequality~(\ref{eq:main}).
We now turn our attention to illustrating the abundance of such trees. Towards this end, we prove that there exist
solutions to inequality~(\ref{eq:main}) with an arbitrary number of vertices of degree~$3$.

Here we use a class of trees called \emph{open quipus} in~\cite{quipu1}. It is shown in~\cite{quipu1} that this class contains examples of infinitely many homeomorphic types to the equation
\[W_2(T)=W(T).\]

\begin{figure}
\begin{center}
\includegraphics[scale=1.25,page=3]{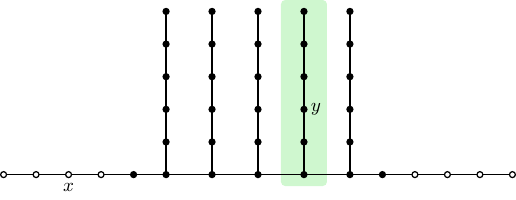}
\end{center}
\caption{The tree $U_{5}$. Deleting the hollow nodes leaves $Q_5$.}
\label{fig:U5}
\end{figure}

For integers $k\ge1$ and $h_1,\ldots,h_k\ge1$, the {open quipu} $Q(k;h_1,\ldots,h_k)$ is constructed on a path $u_0u_1\cdots u_{k+1}$,  by attaching at each node $u_i$ ($1\le i\le k$) a path of $h_i$~vertices.
For an integer $a\ge2$, we denote the open quipu $Q(a;a,\ldots,a)$ by~$Q_a$, and
define the tree $U_a$ to be obtained from $Q_a$ by subdividing each of the edges $u_0u_1$ and $u_a u_{a+1}$ to a path of length~$a$. Figure~\ref{fig:U5} shows~$U_5$.

The next lemma follows from calculations in~\cite{quipu1}.

\begin{lemma}
For all $a\ge2$,
\[W(Q_a)=\frac{2}{3} a^{5}+a^{4}+2 a^{3}+\frac{5}{2} a^{2}+\frac{11}{6} a +1,\]
and
\[D_2(Q_a)=\frac{1}{6} a^{4}+3 a^{3}+\frac{19}{3} a^{2}-\frac{9}{2} a +1.\]
\label{lem:Qa}
\end{lemma}

The next lemma uses these results to estimate $W(U_a)$ and $D_2(U_a)$.

\begin{lemma}
If $a$ is large enough, then
\[W(U_a)=\frac23 a^5 + O(a^4),\]
and
\[D_2(U_a)=\frac16 a^4 + O(a^3).\]
\label{lem:Ua}
\end{lemma}

\begin{figure}
\begin{center}
\includegraphics[scale=1.25,page=4]{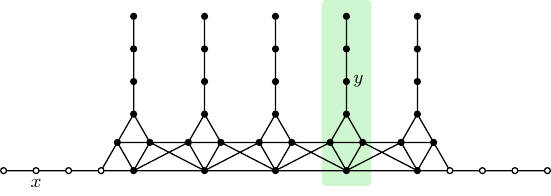}
\end{center}
\caption{The second-order iterated line graph $L^2(U_{5})$. Deleting the hollow nodes leaves $L^2(Q_5)$.}
\label{fig:L2U5}
\end{figure}

\begin{proof}
Using the notation of the definition of $Q_a$, the tree $U_a$ is obtained from $Q_a$ by attaching a path of order $a-1$ at each of the vertices $u_0$ and $u_{a+1}$ (see Figure~\ref{fig:U5}). Then
\begin{equation}
W(U_a)=W(Q_a)+W(P_{3a})-W(P_{a+2})+2\sum\limits_{x=1}^{a-1} \sum\limits_{i=1}^{a} \sum\limits_{y=1}^{a}(x+y+i),
\label{eq:WUa}
\end{equation}
where the term $W(P_{3a})-W(P_{a+2})$ accounts for the increase in the length of the ``spine'' of the tree from $Q_a$ to $U_a$,
and the triple sum adds in the total the distances between any vertex $x\in V(U_a)\setminus V(Q_a)$ and any vertex $y\in V(Q_a)$,
where $i$ indicates the ``arm'' where $y$ resides. Note that $W(P_{3a})-W(P_{a+2})=O(a^3)$, and by direct calculation,
\[2\sum\limits_{x=1}^{a-1} \sum\limits_{i=1}^{a} \sum\limits_{y=1}^{a}(x+y+i)=3a^4-a^3-2a^2=O(a^4).\]
Hence by Lemma~\ref{lem:Qa} we obtain \[W(U_a)=\frac23 a^5+O(a^4).\]

For $D_2(U_a)$, similarly to~(\ref{eq:WUa}), we compute $W_2(U_a)$ by adding to $W_2(Q_a)$ the distances between all unordered pairs $\{x,y\}$ of vertices where $x,y\in V(U_a)\setminus V(Q_a)$, or $x\in V(U_a)\setminus V(Q_a)$ and $y\in V(Q_a)$. See the illustration of $L^2(U_5)$ in Figure~\ref{fig:L2U5} as a reference.
\begin{equation*}
\begin{aligned}
W_2(U_a)=W_2(Q_a)+2W(P_{a-1})&+\sum\limits_{x=1}^{a-1}\sum\limits_{y=1}^{a-1}(x+y+a-1)\\
&+2\sum\limits_{x=1}^{a-1} \sum\limits_{i=1}^{a}\Big[(x+i)+(x+i)+(x+i+1)+ \sum\limits_{y=1}^{a-1}(x+y+i)\Big].
\end{aligned}
\end{equation*}
Direct computation of the nested sums yields $W_2(U_a)=W_2(Q_a)+3a^4+O(a^3)$, which together with Lemma~\ref{lem:Qa} gives
\[D_2(U_a)=\frac16 a^4 + O(a^3).\qedhere\]
\end{proof}

\begin{theorem}
If $a$ is large enough, then $R_2(U_a)<R_2(P_n)$ where $n=|V(U_a)|=a^2+3a$.
\label{thm:Ua}
\end{theorem}

\begin{proof}
Let $n=a^2+3a$. Then $6(n-1)=6a^2+O(n)<7a^2$ for large enough $a$, and $n(n+1)>a^4$ for all~$a\ge1$.
Thus equation~\ref{eq:Pn} gives
\[1-R_2(P_n)=\frac{6(n-1)}{n(n+1)}<\frac7{a^2}.\]
On the other hand, by Lemma~\ref{lem:Ua}, $W(U_a)=\frac23a^5+O(a^4)<a^5$ for large enough~$a$, and
$D_2(U_a)=\frac16a^4+O(a^3)<\frac17a^4$ for large enough~$a$. Therefore,
\[1-R_2(U_a)=\frac{D_2(U_a)}{W(U_a)}>\frac1{7a}.\]
Now if $a\ge50$ is large enough for all previous computations of this proof, we obtain
\[a>49\implies \frac7{a^2}<\frac1{7a}\implies 1-R_2(U_a)>1-R_2(P_n)\implies R_2(U_a)<R_2(P_n).\qedhere\]
\end{proof}

As seen in the proof of Theorem~\ref{thm:star}, if $T=T_{a,b,c}$ is balanced and of order $n$, then
$1-R_2(T)$ and $1-R_2(P_n)$ have the same asymptotic behavior as $n\to\infty$. On the other hand, 
proof of Theorem~\ref{thm:Ua} shows that if $T=U_a$ has order $n$, then as $n\to\infty$,
$1-R_2(P_n)$ approaches $0$ much more quickly than $1-R_2(T)$. Therefore, for orders $n=a^2+3a$ where $a$
is a large enough integer, $R_2(U_a)$ is smaller than the same ratio for the path and the balanced $T_{a,b,c}$ of the same order. 
The construction of $U_a$ can be modified to allow path lengths of $a$ or $a+1$. It is not far-fetched to expect
the resulting tree to behave similarly to $U_a$ in terms of the ratio~$R_2$.

\section{Concluding remarks}

In this work, we give two families of trees which satisfy the inequality~(\ref{eq:main}).
While the trees $T_{a,b,c}$ contains an example of every order $n\ge22$, the quipu-like tree $U_a$
of the same order (if one exists) outperforms them asymptotically.
While these examples constitute infinitely many homeomorphism classes, they all have maximum degree~$3$.
Thus the following problem naturally arises.

\begin{problem}
Is there any tree $T$ with $\Delta(T)\ge4$ that satisfies inequality~(\ref{eq:main})?
\end{problem}

Attempts at higher degrees would naturally start at trees homeomorphic to $K_{1,4}$. It can be shown similarly to
our study of the trees $T_{a,b,c}$ that if $T$ is a ``balanced'' tree homeomorphic to $K_{1,4}$, then $R_2(T)>R_2(P)$, where $P$ is the path of the same order as~$T$.

Let $R_2(n)$ denote the smallest value of $R_2(G)$ among all connected graphs $G$ of order~$n$,
and $R^t_2(n)$ denote the smallest value of $R_2(T)$ among all trees $T$ of order~$n$.
The problems of finding these two extremal values and characterizing the graphs which attain them
remain open. We suggest the following sub-problems as future directions for research.

\begin{problem}
Does $R_2(n)=R^t_2(n)$ hold?
\end{problem}

\begin{problem}
What is the asymptotic behavior of $R_2(n)$ and $R^t_2(n)$? Is it true that they both approach~$1$ as $n\to\infty$?
\end{problem}


\section*{Conflicts of Interest}

Authors declare that they have no conflict of interest.

\def\cprime{$'$}

\end{document}